\documentclass[12pt,a4paper,reqno]{amsart}
\usepackage{graphics,epic}
\usepackage{amsmath,amssymb, amsthm}
\usepackage[all,2cell]{xy}

\newtheorem{theorem}{Theorem}[section]
\newtheorem*{theorem*}{Theorem}

\newtheorem{lemma}[theorem]{Lemma}
\newtheorem{proposition}[theorem]{Proposition}
\newtheorem{corollary}[theorem]{Corollary}

\newtheorem*{conjecture*}{Conjecture}

\newtheorem{remark}[theorem]{Remark}

\renewcommand{\hat}[1]{\widehat{#1}}

\newcommand{\id}{{\rm id}}

\newcommand{\End}{{\rm End}\,}
\newcommand{\Res}{{\rm Res}\,}

\newcommand{\Z}{\mathbb{Z}}

\newcommand{\C}{\mathbb{C}}
\newcommand{\N}{\mathbb{N}}
\def\Res{{\rm Res}}
\def\wt{{\rm wt}}

\def\C{{\mathbb C}}

\def\R{{\mathbb R}}

\def\Z{{\mathbb Z}}

\def\N{{\mathbb N}}

\def\1{{\bf 1}}

\def \End{{\rm End}}

\def \pf{\noindent {\bf Proof: \,}}
\def\theequation{5.\arabic{equation}}

\def \h{\mathfrak{h}}

\def \w{\omega}
\def \g{\mathfrak{g}}

\begin{document}

\title[On the unitary structures of VOSAs]{On the unitary structures of vertex operator superalgebras}
\author{Chunrui Ai}
\address{Chunrui Ai, School of Mathematics and Statistics, Zhengzhou University, Henan 450001, China
}
\email{aicr@zzu.edu.cn}
\author{Xingjun Lin}
\address{Xingjun Lin, Institute of Mathematics, Academia Sinica, Taipei 10617, Taiwan}
\email{linxingjun@math.sinica.edu.tw}
\begin{abstract}
In this paper, the notion of unitary vertex operator superalgebra is introduced. It is proved that the vertex operator superalgebras associated to the unitary highest weight representations for the Neveu-Schwarz Lie superalgebra, Heisenberg superalgebras and  to positive definite integral lattices are unitary vertex operator superalgebras.  The unitary structures are then used to study the structures of vertex operator superalgebras, it is proved that any unitary vertex operator superalgebra is a direct sum of strong CFT type unitary simple vertex operator superalgebras. The classification of unitary vertex operator superalgebras generated by the subspaces with conformal weights less than or equal to $1$ is also considered.
\end{abstract}
\maketitle
\section{Introduction \label{intro}}
\def\theequation{1.\arabic{equation}}
\setcounter{equation}{0}
Unitary structures of vertex operator algebras were introduced in the early days of vertex operator algebras, the unitary structures of the lattice vertex operator algebras and Moonshine vertex operator algebra were used to study the Monster group \cite{B}, \cite{FLM}. Later, based on a symmetric contravariant bilinear form with respect to a Cartan involution in \cite{B} for a vertex  algebra constructed from an even lattice, the notion of invariant bilinear form was introduced and studied in \cite{FHL}, \cite{Li1}.

Recently, it was found in \cite{DLin} that the positive definite Hermitian forms of vertex operator algebras which are invariant with respect to anti-linear involutions of  vertex operator algebras can be used to define unitary vertex operator algebras. And it was proved in \cite{DLin} that the vertex operator algebras associated
to the unitary highest weight representations for the Heisenberg algebras, Virasoro algebra and affine Kac-Moody algebras are unitary vertex operator algebras. Moreover, the unitary structures of these vertex operator algebras are induced from the unitary structures of the highest weight modules for the corresponding Lie algebras. The unitary structures of vertex operator algebras were later used to construct conformal nets from  vertex operator algebras \cite{C}.

In the first part of this paper, the notion of unitary vertex operator superalgebra is introduced, this is  a generalization of the notion of unitary vertex operator algebra.
It is then proved that the vertex operator superalgebras associated to the unitary highest weight representations for the Neveu-Schwarz Lie superalgebra, Heisenberg superalgebras and to positive definite integral lattices are unitary vertex operator superalgebras. Also, the unitary structures of the vertex operator superalgebras associated to the Neveu-Schwarz Lie superalgebra, Heisenberg superalgebras are induced from the unitary structures of the highest weight modules for the corresponding Lie superalgebras. It is also expected that the unitary structures of vertex operator superalgebras can be used to construct superconformal nets (cf. \cite{CKL}).

In the second part of this paper, the structures of unitary vertex operator superalgebras are studied. It is well-known that a finite dimensional associative algebra  over the field of complex numbers is semisimple if and only if the trace form of the algebra is nondegenerate and that any semisimple associative algebra over the field of complex numbers is a direct sum of simple algebras. Since vertex operator algebras are analogues of associative algebras and unitary vertex operator superalgebras have positive definite invariant  Hermitian forms, it is natural to expect that any unitary vertex operator superalgebra is a direct sum of unitary simple vertex operator superalgebras. In this paper, it is shown that any unitary vertex operator superalgebra is a direct sum of strong CFT type unitary simple vertex operator superalgebras. And it is found that the strong CFT type properties of vertex operator superalgebras, which have played important roles in the study of vertex operator superalgebras (see \cite{Bu}, \cite{DM1}, \cite{DM2}, \cite{Hu}), are closely related to the unitary structures of vertex operator superalgebras.

As an application of the structure theorem of unitary vertex operator superalgebras, the classification of unitary vertex operator superalgebras generated by the subspaces with conformal weights less than or equal to $1$ is considered. The key point is to classify the simple vertex operator superalgebras  in these classes. The simple vertex operator superalgebras generated by the subspaces with conformal weight $\frac{1}{2}$ have been classified in \cite{X}, it was proved that any simple vertex operator superalgebra generated by the subspace with conformal weight $\frac{1}{2}$ is isomorphic to the vertex operator superalgebra associated to the highest weight representation for some Heisenberg superalgebra. So the key point is to classify unitary simple vertex operator algebras generated by the subspaces with conformal weight $1$. Simple vertex operator algebras generated by the subspaces with conformal weight $1$ have been studied in \cite{DM1}, \cite{DM2}, \cite{M} under the condition that the vertex operator algebras are $C_2$-cofinite or regular, and their results have played important roles in the classification of rational and $C_2$-cofinite vertex operator algebras. But their results depend on the condition that the vertex operator algebras are $C_2$-cofinite or regular, and there are unitary vertex operator algebras which are neither $C_2$-cofinite nor regular, so we still need to study the structures of unitary simple vertex operator algebras generated by the subspaces with conformal weight $1$. In this paper, it is shown that each of these unitary vertex operator superalgebras is isomorphic to the tensor product of some vertex operator superalgebras associated to unitary highest weight representations for the affine Kac-Moody algebras and the Heisenberg algebras.

The rest of this paper is organized as follows: In Section 2, we prove some basic facts about unitary vertex operator superalgebras and then show that the vertex operator superalgebras associated to the unitary highest weight representations for the Neveu-Schwarz Lie superalgebra, Heisenberg superalgebras and to positive definite integral lattices are unitary vertex operator superalgebras. In Section 3, the structures of unitary vertex operator superalgebras are studied. In Section 4, the classification of unitary vertex operator superalgebras generated by the subspaces with conformal weights less than or equal to $1$ is considered.
\section{Unitary vertex operator superalgebras}
\def\theequation{2.\arabic{equation}}
\setcounter{equation}{0}
\subsection{Preliminaries} In this subsection, we shall recall from \cite{DL}, \cite{FFR}, \cite{KW}, \cite{Li} and \cite{X} some facts about vertex operator superalgebras. First, we recall some notations, let $V=V_{\bar 0}\oplus V_{\bar 1}$ be any ${\Bbb Z}_2$-graded vector space, the element in $V_{\bar 0} $ (resp. $V_{\bar 1}$) is called  {\em even} (resp. {\em odd}).
We then define $[v]=i$ for any $v\in V_{\bar i}$ with  $i=0,1$.  A {\em vertex operator superalgebra} is a quadruple $(V,Y,\1,\w),$ where $V=V_{\bar 0}\oplus V_{\bar 1}$ is a ${\Bbb Z}_2$-graded vector space, ${\bf 1}$ is the {\em vacuum vector}
of $V$, $\w$ is the {\em conformal vector} of $V,$   and $Y$ is a linear map
\begin{align*}
  V &\to (\End\,V)[[z,z^{-1}]] ,\\
 v&\mapsto Y(v,z)=\sum_{n\in{\Z}}v_nz^{-n-1}\ \ \ \  (v_n\in
\End\,V)
\end{align*}
satisfying the following axioms:  \\
(i) For any $u,v\in V,$ $u_nv=0$ for sufficiently large $n$;\\
(ii) $Y({\bf 1},z)=\id_{V}$;\\
(iii) $Y(v,z){\bf 1}=v+\sum_{n\geq 2}v_{-n}{\bf 1}z^{n-1},$ for any $v\in V$;\\
(iv) The component operators of  $Y(\w,z)=\sum_{n\in\Z}L(n)z^{-n-2}$ satisfy the Virasoro algebra
relation with {\em central charge} $c\in \C:$
\begin{align*}
[L(m),L(n)]=(m-n)L(m+n)+\frac{1}{12}(m^3-m)\delta_{m+n,0}c,
\end{align*}
and
\begin{align*}
\frac{d}{dz}Y(v,z)=Y(L(-1)v,z),\ \ \text{ for any } v\in V;
\end{align*}
(v) $V$ is $\frac{1}{2}\Z$-graded such that $V=\oplus_{n\in \frac{1}{2}\Z}V_n$, $V_{\bar 0}=\oplus_{n\in \Z}V_n$, $V_{\bar 1}=\oplus_{n\in \frac{1}{2}+\Z}V_n$, $L(0)|_{V_n}=n$, $\dim(V_n)<\infty$ and $V_n=0$ for sufficiently small $n$;\\
(vi) The {\em Jacobi identity} for ${\Bbb Z}_2$-homogeneous $u,v\in V$ holds,
\begin{align*}
\begin{array}{c}
\displaystyle{z^{-1}_0\delta\left(\frac{z_1-z_2}{z_0}\right)
Y(u,z_1)Y(v,z_2)-(-1)^{[u][v]} z^{-1}_0\delta\left(\frac{z_2-z_1}{-z_0}\right)
Y(v,z_2)Y(u,z_1)}\\
\displaystyle{=z_2^{-1}\delta
\left(\frac{z_1-z_0}{z_2}\right)
Y(Y(u,z_0)v,z_2)}.
\end{array}
\end{align*}
This completes the definition of a vertex operator superalgebra and we will denote the vertex operator superalgebra briefly by $V$. An {\em anti-linear automorphism} $\phi$ of $V$ is
an anti-linear isomorphism (as anti-linear morphism) $\phi:V\to V$ such
that $\phi(\bf1)=\bf1, \phi(\omega)=\omega$ 
and
$\phi(u_nv)=\phi(u)_n\phi(v)$ for any $ u, v\in V$ and $n\in
\mathbb{Z}$. An anti-linear automorphism $\phi$ of $V$ is called an {\em anti-linear involution} if the order of $\phi$ is $2$.
We then define a {\em unitary} vertex operator superalgebra to be a pair $(V,\phi)$, where $V$ is a vertex operator superalgebra and $\phi: V\to V$ is an anti-linear involution of $V$,  such that there exists a positive definite Hermitian form $(,): V\times V\to \mathbb{C}$ on $V$ and the following invariant property 
 $$(Y(e^{zL(1)}(-1)^{L(0)+2L(0)^2}z^{-2L(0)}a, z^{-1})u, v)=(u,
Y(\phi(a), z)v),$$holds for any $a, u, v\in V$.
\begin{remark}
Note that if $V$ is a vertex operator algebra, the definition of unitary vertex operator superalgebra coincides with the definition of unitary vertex operator algebra in \cite{DLin}.
\end{remark}

Similarly, we can define the unitary module for a unitary vertex operator superalgebra (cf. \cite{DLin}). The following results are immediate consequences of the definitions of unitary modules and unitary vertex operator superalgebras.
\begin{proposition}
  Let $V$ be a unitary vertex operator superalgebra $V$. Then \\
  (1) For any homogenous vectors $u,v \in V$ such that $\wt u\neq \wt v$, we have $(u, v)=0$.\\
  (2) Any unitary $V$-module $M$  is completely reducible, that is, $M \cong \oplus_{i\in I}W^i$, where $W^i$, $i\in I$, are irreducible submodules of $M$. Moreover, $(W^i, W^j)_M=0$ for  $i\neq j$, where $(,)_M$ is the positive definite Hermitian form on $M$.
\end{proposition}
\pf We only prove the second result. Let $M$ be a unitary $V$-module and $(,)_M$ be the positive definite Hermitian form on $M$. Consider the set $S$ consisting of submodules of $M$ such that each element $W\in S$ satisfies the property $W \cong \oplus_{i\in I}W^i$,  and $(W^i, W^j)_M=0$ for $ i\neq j$, where $W^i, i\in I$, are irreducible submodules of $M$.
Note that there exists a maximal object $W^0$ in $S$. We now prove that $W^0=M$. Otherwise, consider the orthogonal complement $(W^0)^\perp$ of $W^0$ with respect to $(,)_M$. It follows from the definition of unitary modules that $(W^0)^\perp$ is also a submodule of $M$. Let $n_0$ be a number such that $(W^0)^\perp_{n_0}\neq 0$ but $(W^0)^\perp_{n}= 0$ if $Re~ n< Re~ n_0$, where $(W^0)^\perp_{n}=\{w\in (W^0)^\perp|L(0)w=nw\}$. Then we know that $(W^0)^\perp_{n_0}$ is a module for the Zhu's algebra $A(V)$ of $V$ (cf. \cite{KWa}). Since $\dim (W^0)^\perp_{n_0}<\infty$, there exists an irreducible $A(V)$-submodule $Z$ of $(W^0)^\perp_{n_0}$, we then prove that the submodule $\langle Z\rangle$  of $(W^0)^\perp$ generated by $Z$ is an irreducible module for $V$. Note that $\langle Z\rangle\cap (W^0)^\perp_{n_0}=Z$ and that any proper submodule of $\langle Z\rangle$ has zero intersection with $Z$, it follows immediately from the invariant property of $(,)_M$ that $\langle Z\rangle$ is irreducible. 
As a result, $W^0\oplus \langle Z\rangle$ is an object in $S$, this is a contradiction. The proof is complete. \qed

\vskip.25cm
 Recall that for a vertex operator superalgebra $V$ and a $V$-module $M=\oplus_{\lambda\in \C}M_{\lambda}$,  the {\em  contragredient module} $M'$ of $M$ is defined to be $(M', Y')$, where $M'=\oplus_{\lambda\in \C} (M_\lambda)'$ is the restricted dual of $M$ and $Y'$ is the linear map defined by
  \begin{align*}
  \langle Y'(a, z)w', w\rangle = \langle w', Y(e^{zL(1)}(-1)^{L(0)+2L(0)^2}z^{-2L(0)}a, z^{-1})w\rangle,
  \end{align*}for any $a\in V, w\in M, w'\in M'$ (cf. \cite{FHL}, \cite{Y}). Recall also that a $V$-module $M$ is called {\em self-dual} if $M$ is isomorphic to $M'$, then we have
\begin{proposition}\label{self-dual}
Let $V$ be a unitary vertex operator superalgebra. Then $V$ is self-dual.
\end{proposition}
\pf  It is good enough to prove that $V$ has a nondegenerate invariant bilinear form, that is, a nondegenerate bilinear form $\langle ,\rangle$ such that the following invariant property 
\begin{align*}
\langle Y(u, z)v, w\rangle=\langle v, Y(e^{zL(1)}(-1)^{L(0)+2L(0)^2}z^{-2L(0)}u, z^{-1})w\rangle,
\end{align*}
holds for any $u,v, w\in V$. For any $u,v\in V$, we define $\langle u, v\rangle =(u, \phi(v))$. It is clear that $\langle, \rangle$ is a nondegenerate bilinear form. We then prove that $\langle, \rangle$ is invariant. By definition, we have
\begin{align*}
\langle Y(u, z)v, w\rangle
& =(Y(u, z)v, \phi (w))\\
& =(v, Y(e^{zL(1)}(-1)^{L(0)+2L(0)^2}z^{-2L(0)}\phi(u), z^{-1})\phi (w))\\
& =(v, \phi(Y(e^{zL(1)}(-1)^{L(0)+2L(0)^2}z^{-2L(0)}u, z^{-1})w))\\
& =\langle v, Y(e^{zL(1)}(-1)^{L(0)+2L(0)^2}z^{-2L(0)}u, z^{-1})w\rangle,
\end{align*}
as desired.\qed

The following results concern constructing unitary vertex operator superalgebras from known unitary vertex operator superalgebras, the proofs of these results are similar to those of Propositions 2.6, 2.9 in \cite{DLin}.
\begin{proposition}\label{sub}
(1) Let $(V, \phi)$ be a unitary vertex operator superalgebra and $U$ be a subalgebra of $V$ such that the Virasoro element of $U$ is the same as that of $V$ and $\phi(U)=U$. Then $(U, \phi|_U)$ is a unitary vertex operator superalgebra.\\
(2) Let $(V^1, \phi_1),..., (V^p, \phi_p)$ be unitary vertex operator superalgebras. Then there exists an anti-linear involution $\phi$ of $V^1\otimes\dots\otimes V^p $ such that  $(V^1\otimes\dots\otimes V^p, \phi)$ is a unitary vertex operator superalgebra.
\end{proposition}

We also need the following result which is useful to prove the unitarity of vertex operator superalgebra.
\begin{proposition}\label{basic}
Let $V$ be a vertex operator superalgebra equipped with a positive
definite Hermitian form $(, ): V\times V\to \mathbb{C}$ and $\phi$ be an anti-linear involution of $V$. Assume that $V$ is
generated by the subset $S\subset V$, i.e.,
$$V=span\{u^1_{n_1}\cdots u^k_{n_k}{\bf 1}| k\in \N, u^1,\cdots,
u^k\in S\}$$ and that the invariant property
\begin{equation*} (Y(e^{zL(1)}(-1)^{L(0)+2L(0)^2}z^{-2L(0)}a, z^{-1})u,
v)=(u, Y(\phi(a), z)v)
 \end{equation*}
 holds for any $a\in S, u, v\in V.$ Then $(V, \phi)$ is a unitary vertex operator superalgebra.
\end{proposition}
\pf Note that for a vertex operator superalgebra $V$, we have the following identity which was essentially proved in  Theorem 5.2.1 of \cite{FHL},
\begin{align*}
&-(-1)^{[a][b]}z _0^{-1}\delta\left(\frac{z_2-z_1}{-z_0}\right)Y(e^{z_1L(1)}(-1)^{L(0)+2L(0)^2}z_1^{-2L(0)}a, z_1^{-1})\\
&~~\cdot Y(e^{z_2L(1)}(-1)^{L(0)+2L(0)^2}z_2^{-2L(0)}b, z_2^{-1})\\
&+z_0^{-1}\delta\left(\frac{z_1-z_2}{z_0}\right)Y(e^{z_2L(1)}(-1)^{L(0)+2L(0)^2}z_2^{-2L(0)}b, z_2^{-1})\\
&~~\cdot Y(e^{z_1L(1)}(-1)^{L(0)+2L(0)^2}z_1^{-2L(0)}a, z_1^{-1})\\
&\ \ =z_1^{-1}\delta\left(\frac{z_2+z_0}{z_1}\right)Y(e^{z_2L(1)}(-1)^{L(0)+2L(0)^2}z_2^{-2L(0)}Y(a, z_0)b, z_2^{-1}),
\end{align*}
for any $a, b\in V$. Then the proposition follows from the similar argument in the proof of  Proposition 2.11 in \cite{DLin}. \qed
\subsection{Unitary Neveu-Schwarz vertex operator superalgebras}
In this subsection, we shall show that the Neveu-Schwarz vertex operator superalgebras associated to unitary highest weight modules for the Neveu-Schwarz Lie superalgebra are unitary vertex operator superalgebras. We first recall from \cite{KWa}, \cite{Li} some facts about the Neveu-Schwarz vertex operator superalgebras. Let $NS$ be the Neveu-Schwarz Lie superalgebra
\begin{align*}
NS=\oplus_{m\in \Z}\C L(m)\oplus\oplus_{n\in \Z}\C G(n+\frac{1}{2})\oplus \C C,
\end{align*}
with the following communication relations:
\begin{align*}
&[L(m), L(n)]=(m-n)L(m+n)+\frac{m^3-m}{12}\delta_{m+n, 0}C,\\
&[L(m), G(n+\frac{1}{2})]=(\frac{m}{2}-n-\frac{1}{2})G(m+n+\frac{1}{2}),\\
&[G(m+\frac{1}{2}), G(n+\frac{1}{2})]_+=2L(m+n)+\frac{1}{3}m(m+1)\delta_{m+n, 0}C,\\
&[L(m), C]=0,~[G(n+\frac{1}{2}), C]=0.
\end{align*}
Set
\begin{align*}
&NS_{\pm}=\oplus_{n\in \Z_+}\C L(\pm n)\oplus\oplus_{m\in \Z_+}\C G(\pm m\mp\frac{1}{2}),\\
&NS_0=\C L(0)\oplus \C C.
\end{align*}
It is clear that $NS_+\oplus NS_0$ is a subalgebra of $NS$. Given any complex numbers $c$ and $h$,  we then consider the Verma module
\begin{align*}
M(c, h)=U(NS)/J,
\end{align*}
 where $U(\g)$ denotes the universal enveloping algebra of a Lie superalgebra $\g$ and  $J$ is the left ideal of $U(NS)$ generated by $NS_+$, $L(0)-h$, $C-c$.
 It is well-known that the Verma module has a unique maximal proper submodule $J(c, h)$, we then let $L(c, h)$ be the irreducible $NS$-module $M(c, h)/J(c, h)$. It was proved in \cite{KWa} that $L(c, 0)$ has a vertex operator superalgebra structure such that $\1=1$ and $\w=L(-2)1$.

 We next recall from \cite{KW} some facts about the unitary modules for the Neveu-Schwarz Lie superalgebra $NS$. Let $M$ be a highest weight module for $NS$, $M$ is called a {\em unitary} module for $NS$ if there exists a positive definite Hermitian form $(,)$ on $M$ such that
 \begin{align*}
 (L(n)u,v)=(u, L(-n)v),~(G(n+\frac{1}{2})u, v)=(u,G(-n-\frac{1}{2})v ),
 \end{align*}
 for any $n\in \Z$ and $u,v \in M$ (cf. \cite{KT}, \cite{KW}). It was proved in \cite{GKO}, \cite{KW} that $L(c, h)$ is a unitary module for $NS$ if $c\geq \frac{3}{2}$ and $h\geq 0$ or
 \begin{align*}
 &c=c_m=\frac{3}{2}(1-\frac{8}{(m+2)(m+4)}),\\
 &h=h^m_{r, s}=\frac{((m+4)r-(m+2)s)^2-4}{8(m+2)(m+4)},
 \end{align*}
 where $m, r,s\in \Z_+,~ 1\leq s\leq r\leq m+1$ and $r-s\in 2\Z$ if $r\neq 0$.

 We now begin to construct unitary structures on the vertex operator superalgebras associated to unitary modules for $NS$. We begin with the definition of the anti-linear involution. Let $T(NS)$ be the tensor algebra of $NS$ and $\Phi: T(NS)\to T(NS)$ be the anti-linear map of $T(NS)$ such that $\Phi(x^1\otimes\cdots \otimes x^n)=x^1\otimes\cdots \otimes x^n$ for any $x^1,\cdots , x^n\in \{L(n), G(m+\frac{1}{2})|n,m\in \Z\}$.
   It is clear that $\Phi$ is an anti-linear involution of $T(NS)$. Note that by the definition of $NS$,
   we have\begin{align*}
  &\Phi(x^1\otimes x^2-(-1)^{[x^1][x^2]}x^2\otimes x^1-[x^1, x^2])\\
  &\ \ =\Phi(x^1)\otimes \Phi(x^2)-(-1)^{[x^1][x^2]}\Phi(x^2)\otimes\Phi( x^1)-\Phi([x^1, x^2])\\
  &\ \ =\Phi(x^1)\otimes \Phi(x^2)-(-1)^{[x^1][x^2]}\Phi(x^2)\otimes\Phi( x^1)-[\Phi(x^1), \Phi(x^2)],
  \end{align*}
   for any $x^1, x^2\in NS$. Thus, $\Phi$ induces an anti-linear involution of $U(NS)$, which is still denoted by $\Phi$. Furthermore, if $c\in \R$, we have $\Phi(J)\subset J$. Hence, $\Phi$ induces an anti-linear involution of $M(c, 0)$, which is also denoted by $\Phi$. Finally, note that $\Phi(J(c, 0))\subset J(c, 0)$, we then have an anti-linear involution $\phi$ of $L(c, 0)$ induced from $\Phi$. We are now in a position to prove that $L(c,0)$ has a unitary vertex operator superalgebra structure.
 \begin{theorem}
 Let $c\in \R$ be a positive real number such that $c\geq \frac{3}{2}$ or $c=c_m$ for some positive integer $m$. Then $(L(c, 0), \phi)$ is a unitary vertex operator superalgebra.
 \end{theorem}
 \pf Note that if $c\in \R$ is a positive real number such that $c\geq \frac{3}{2}$ or $c=c_m$ for some positive integer $m$,  $L(c, 0)$ is a unitary module for $NS$. Let $(,)$ be the positive definite Hermitian form on $L(c, 0)$. We next prove that $(L(c, 0), \phi)$ equipped with the positive definite Hermitian form $(,)$ is a unitary vertex operator superalgebra. Since $L(c, 0)$ is generated by $L(-2)1$ and $G(-\frac{3}{2})1$, by Proposition \ref{basic} it is good enough to prove that the invariant property
 \begin{equation*} (Y(e^{zL(1)}(-1)^{L(0)+2L(0)^2}z^{-2L(0)}a, z^{-1})u,
v)=(u, Y(\phi(a), z)v)
 \end{equation*}
 holds for $a=L(-2)1\text{ or }G(-\frac{3}{2})1$ and any $ u, v\in L(c, 0).$ When $a=L(-2)1$, the invariant property has been proved in \cite{DLin}. We then assume that $a=G(-\frac{3}{2})1$, by a direct computation,
 \begin{align*}
 (Y(e^{zL(1)}(-1)^{L(0)+2L(0)^2}z^{-2L(0)}G(-\frac{3}{2})1, z^{-1})u, v)&=z^{-3}(Y(G(-\frac{3}{2})1, z^{-1})u, v)\\
 &=z^{-3}\sum_{n\in \Z}z^{n+2}(G(n+\frac{1}{2})u, v)\\
 &=\sum_{n\in \Z}z^{n-1}(u,G(-n-\frac{1}{2})v)\\
 &=(u, Y(G(-\frac{3}{2})1, z)v)\\
 &=(u, Y(\phi(G(-\frac{3}{2})1), z)v),
 \end{align*}
 as desired.\qed
 \subsection{Unitary  vertex operator superalgebras associated to the Heisenberg superalgebras}\label{superhe}
 In this subsection, we shall show that vertex operator superalgebras associated to the unitary highest weight representations for the Heisenberg superalgebras are unitary vertex operator superalgebras. First, we recall from \cite{KWa}, \cite{Li}, \cite{T}  and  \cite{X} some facts about highest weight modules for the Heisenberg superalgebras. Let $\h$ be a $n$-dimensional vector space with a nondegenerate symmetric bilinear form $\langle,\rangle$. Consider the infinite dimensional Lie superalgebra
 \begin{align*}
 \hat \h=\h\otimes \C[t, t^{-1}]t^{\frac{1}{2}}\oplus \C C,
 \end{align*}
 with $\Z_2$-gradation $\hat \h_{\bar 0}=\C C$, $\hat \h_{\bar 1}=\h\otimes \C[t, t^{-1}]t^{\frac{1}{2}}$ and the communication relations:
 \begin{align*}
 [u(m), v(n)]_+=\delta_{m+n, 0}\langle u, v\rangle C,~~[u(m), C]=0,
 \end{align*}
 for any $u, v\in \h$, $m, n\in \frac{1}{2}+\Z$, where $u(m)=u\otimes t^m$. Set $\hat \h_+=t^{\frac{1}{2}}\C[t]\otimes \h$ and $\hat \h_-=t^{-\frac{1}{2}}\C[t^{-1}]\otimes \h$, it is clear that $\hat \h_+$ and $\hat \h_-$ are subalgebras of $\hat \h$. For any complex number $c$, we then consider the Verma module
 \begin{align*}
 M_{\hat\h}(c, 0)=U(\hat \h)/J,
 \end{align*}
  where $J$ is the left ideal of $U(\hat \h)$ generated by $\hat \h_+$ and $C-c$. It is well-known that $M_{\hat\h}(c, 0)$ is a highest weight module for $\hat \h$ and that $M_{\hat\h}(c, 0)$ is an irreducible module for $\hat \h$ (cf. \cite{Li}). We then let $u^1, \cdots, u^n$ be an orthonormal basis of $\h$ with respect to $\langle,\rangle$ and set $\1=1$, $\w=\frac{1}{2}\sum_{i=1}^n u^i({-\frac{3}{2}})u^i({-\frac{1}{2}})\1$, it is known that $M_\h(c, 0)$ has a vertex operator superalgebra structure such that $\1$ and $\w$ are the vacuum vector and conformal vector, respectively (cf. \cite{KWa}, \cite{Li}, \cite{X}). Moreover, it was proved in \cite{Li} that the vertex operator superalgebra $M_{\hat\h}(c, 0)$ is isomorphic to the vertex operator superalgebra $M_{\hat\h}(1, 0)$ for any nonzero complex number $c$.

 We now begin to construct a unitary structure on $M_{\hat\h}(1, 0)$. We first assume that $\h$ is one-dimensional. Let $\alpha$ be a vector in $\h$ such that $\langle \alpha, \alpha\rangle =1$, then $\h=\C \alpha$. By the PBW Theorem, we know that $M_{\hat\h}(1, 0)$ has a basis consisting of $\alpha(-n_1)\cdots \alpha(-n_k)\1$, $n_1>\cdots>n_k\geq \frac{1}{2}$. It is known that there exists a positive definite Hermitian form $(, )$ on $M_{\hat\h}(1, 0)$ such that $\alpha(-n_1)\cdots \alpha(-n_k)\1$, $n_1>\cdots>n_k\geq \frac{1}{2}$ is an orthonormal basis and that
  \begin{align*}
  ( \alpha(m)u, v )=( u, \alpha(-m)v)
  \end{align*}
 for any $u, v\in M_{\hat\h}(1, 0)$ (cf. \cite{KR}).

 We next construct an anti-linear involution of $M_{\hat\h}(1, 0)$. Let $T(\hat \h)$ be the tensor algebra of $\hat \h$. Define an anti-linear map $\Phi: T(\hat\h)\to T(\hat\h)$ by  $\Phi(\alpha(n_1)\otimes \cdots\otimes \alpha(n_k))=(-1)^k\alpha(n_1)\otimes \cdots\otimes \alpha(n_k)$ and $\Phi(C)=C$. It is clear that $\Phi$ is an anti-linear involution of $T(\hat \h)$. Note that
   \begin{align*}
  &\Phi(x^1\otimes x^2-(-1)^{[x^1][x^2]}x^2\otimes x^1-[x^1, x^2])\\
  &\ \ =\Phi(x^1)\otimes \Phi(x^2)-(-1)^{[x^1][x^2]}\Phi(x^2)\otimes\Phi( x^1)-\Phi([x^1, x^2])\\
  &\ \ =\Phi(x^1)\otimes \Phi(x^2)-(-1)^{[x^1][x^2]}\Phi(x^2)\otimes\Phi( x^1)-[\Phi(x^1), \Phi(x^2)],
  \end{align*}
   for any $x^1, x^2\in \hat\h$. It follows that $\Phi$ induces an anti-linear involution $\phi$ of $M_{\hat\h}(1, 0)$.

 We now in a position to provide a unitary vertex operator superalgebra structure on $M_{\hat\h}(1, 0)$.
  \begin{theorem}
 Let $\h$ be a $n$-dimensional vector space with a nondegenerate symmetric bilinear form $\langle, \rangle$. Then $M_{\hat\h}(1, 0)$ has a unitary vertex operator superalgebra structure.
 \end{theorem}
 \pf Let $u^1, \cdots, u^n$ be an orthonormal basis of $\h$ with respect to $\langle,\rangle$. Then we know that the vertex operator superalgebra $M_{\hat\h}(1, 0)$ is isomorphic to the vertex operator superalgebra $M_{\hat{\C u^1}}(1, 0)\otimes \cdots\otimes M_{\hat{\C u^n}}(1, 0)$.
 By Proposition \ref{sub}, it is good enough to prove that $M_{\hat{\C u^i}}(1, 0)$ has a unitary vertex operator superalgebra structure for each $i$. We now prove that $(M_{\hat{\C \alpha}}(1, 0), \phi)$ is a unitary vertex operator superalgebra. Since $M_{\hat{\C \alpha}}(1, 0)$ is generated by $\alpha(-\frac{1}{2})\1$,   by Proposition \ref{basic} it is good enough to prove that the invariant property
 \begin{equation*} (Y(e^{zL(1)}(-1)^{L(0)+2L(0)^2}z^{-2L(0)}a, z^{-1})u,
v)=(u, Y(\phi(a), z)v)
 \end{equation*}
 holds for $a=\alpha(-\frac{1}{2})\1$ and any $ u, v\in M_{\hat{\C \alpha}}(1, 0).$ By a direct computation, \begin{align*}
 &(Y(e^{zL(1)}z^{-2L(0)}(-1)^{L(0)+2L(0)^2}\alpha(-\frac{1}{2})\1, z^{-1})u, v)\\
 &\ \ =-z^{-1}(Y(\alpha(-\frac{1}{2})\1, z^{-1})u, v)\\
 &\ \ =-z^{-1}\sum_{m\in \Z}(\alpha(m+\frac{1}{2})u, v)z^{m+1}\\
 &\ \ =-z^{-1}\sum_{m\in \Z}(u, \alpha(-m-\frac{1}{2})v)z^{m+1}\\
 &\ \ =-\sum_{m\in \Z}(u, \alpha(-(m+1)+\frac{1}{2})v)z^{m}\\
 &\ \ =(u, Y(-\alpha(-\frac{1}{2})\1, z)v)\\
 &\ \ =(u, Y(\phi(\alpha(-\frac{1}{2})\1), z)v),
 \end{align*}
 as desired. \qed

\subsection{Unitary lattice vertex operator superalgebras}In this subsection we shall prove that the lattice vertex operator superalgebras associated to positive definite integral lattices are unitary vertex operator superalgebras.
First, we recall from \cite{FLM}, \cite{K1} some facts about the lattice vertex operator
superalgebras. Let $L$ be a positive definite integral lattice and $(,)$ be the associated positive definite bilinear form. 
We then
consider the central extension $\hat{L}$ of $L$ by the cyclic
group $\langle \kappa\rangle$ of order $2$:
$$1\to \langle \kappa\rangle\to \hat{L}\to L\to 1,$$ such that the commutator map
$c(\alpha, \beta)=\kappa^{(\alpha, \beta)+(\alpha, \alpha)(\beta, \beta)}$ for any $\alpha, \beta\in L$. Let $e: L\to \hat{L}$
be a section such that $e_0 = 1$ and $\epsilon_0: L\times L\to \langle \kappa\rangle$ be the corresponding 2-cocycle. Then $\epsilon_0(\alpha,
\beta)\epsilon_0(\beta, \alpha)=\kappa^{(\alpha, \beta)+(\alpha, \alpha)(\beta, \beta)}$  and
$e_{\alpha}e_{\beta}=\epsilon_0(\alpha, \beta)e_{\alpha+\beta}$ for
$ \alpha, \beta \in L$.

Next, we consider the induced
$\hat{L}$-module:$$\C\{L\}=\C[\hat{L}]\otimes_{\langle \kappa\rangle}\C\cong
\C[L]\ \ (\text{linearly}),$$ where $\C[L]$ denotes the group algebra of $L$ and
$\kappa$ acts on $\C$ as multiplication by $-1$.  Then $\C[L]$ becomes an
$\hat{L}$-module such that $e_{\alpha}\cdot
e^{\beta}=\epsilon(\alpha, \beta)e^{\alpha+\beta}$ and $\kappa\cdot
e^{\beta}=-e^{\beta}$, where $\epsilon(\alpha, \beta)$ is defined as $\nu\circ\epsilon_0$ and $\nu$ is the an isomorphism from $\langle \kappa\rangle$ to $\langle \pm 1\rangle$ such that $\nu(\kappa)=-1$. We also define an action $h(0)$ on
$\C[L]$ by $h(0)\cdot e^{\alpha}=(h, \alpha)e^{\alpha}$ for $h
\in \h, \alpha\in {L}$ and an action $z^h$ on $\C[L]$ by $z^{h}\cdot e^{\alpha}=z^{(h,
\alpha)}e^{\alpha}$.

We proceed to construct the vector space of the lattice vertex operator superalgebra. Set $\h=\C\otimes_{\Z}L$,  and consider the Heisenberg algebra $\hat{\h}=\h\otimes \C[t, t^{-1}]\oplus \C K$ with the communication relations: For any $u, v\in \h$,
\begin{align*}
&[u(m), v(n)]=m\langle u, v\rangle\delta_{m+n, 0}K,\\
&[K, x]=0, \text{ for any } x\in \hat{\h},
\end{align*}
where $u(n)=u\otimes t^n$, $u\in \h$. 
Let $M(1)$ be the
Heisenberg vertex operator algebra associated to $\hat\h$. The vector space of the lattice vertex operator superalgebra is  defined to be
$$V_L=M(1)\otimes_{\C}\C\{L\}\cong M(1)\otimes_{\C}\C[L]\ \ (\text{linearly}).$$Then $\hat{L}$, $h(n)(n\neq 0)$, $h(0)$ and $z^{h(0)}$ act naturally on $V_L$ by acting on either $M(1)$ or $\C[L]$ as indicated above. It was proved in  \cite{FLM}, \cite{K1} that
$V_L$ has a vertex operator superalgebra structure such that
\begin{align*}
&Y(h(-1)1, z)=h(z)=\sum_{n\in \Z}h(n)z^{-n-1}\ \ (h\in \h),\\
&Y(e^{\alpha}, z)=E^-(-\alpha, z)E^+(-\alpha, z)e_{\alpha}z^{\alpha},
\end{align*}
where
\begin{align*}
E^-(\alpha, z)=\exp\left(\sum_{n<0}\frac{\alpha(n)}{n}z^{-n}\right),\ \ \
E^+(\alpha, z)=\exp\left(\sum_{n>0}\frac{\alpha(n)}{n}z^{-n}\right).
\end{align*}

We now begin to construct a unitary vertex operator superalgebra structure on $V_L$. First, we define a positive definite Hermitian form on $V_{L}$. Recall from \cite{DLin} that there is
a unique positive definite Hermitian form $(,)$ on $M(1)$ such
that
$$({\bf 1}, {\bf 1})=1,\ \ \ \ (h(n)u, v)=(u,
h(-n)v)$$
for  any $h\in \h, u, v\in M(1).$ Also, there is a positive definite Hermitian form
$(,):\C[L]\times\C[L]\to \C$
on
$\C[L]$ determined by   $(e^{\alpha},e^{\beta})=0$ if
$\alpha\neq\beta$ and $(e^{\alpha},e^{\beta})=1$ if
$\alpha=\beta$ (see \cite{FLM}).
We then define a positive definite Hermitian form $(,)$ on $V_{L}$
as follows: For any  $u, v\in M(1)$ and $e^{\alpha}, e^{\beta}\in \C[L]$, $$(u\otimes e^{\alpha}, v\otimes e^{\beta})=(u,
v)(e^{\alpha}, e^{\beta}).$$

\begin{lemma}\label{lattv}
Let $(,)$ be the positive definite Hermitian form on $V_{L}$ defined above.
Then we have: For any $\alpha \in L$ and  $ w_1, w_2\in
V_{L}$,
\begin{align*}
(e_{\alpha}w_1, w_2)=(w_1,(-1)^{\frac{(\alpha, \alpha)+(\alpha,\alpha)^2}{2}} e_{-\alpha}w_2),\ \ \ \ (z^{\alpha}w_1, w_2)=(w_1, z^{\alpha}w_2).
\end{align*}
\end{lemma}
\pf The second identity is obvious. The first identity follows immediately from the facts that $(e_{\alpha} w_1, e_{\alpha}w_2)=(w_1, w_2)$ for any $w_1, w_2\in V_{L}$ and that $e_\alpha e_{-\alpha}=(-1)^{\frac{(\alpha, \alpha)+(\alpha, \alpha)^2}{2}}$ for any $\alpha\in L$ (see \cite{K1}).\qed

We next construct an anti-linear involution of $V_L$. Let $\phi: V_L\to V_L$ be an anti-linear map  determined by:
\begin{align*}
\phi: V_L &\to V_L\\
\alpha_1(-n_1)\cdots\alpha_k(-n_k)\otimes e^{\alpha}&\mapsto (-1)^k\alpha_1(-n_1)\cdots\alpha_k(-n_k)\otimes
e^{-\alpha}.
\end{align*}
Note that we have $$\phi\alpha(n)\phi^{-1}=-\alpha(n),\ \ \ \ \phi Y(e^{\alpha}, z)\phi^{-1}=Y(e^{-\alpha}, z),$$
 for any $\alpha\in L$, $n\in \Z.$
By a similar argument as that of Lemma 4.1 in \cite{DLin}  we can show that $\phi$ is an anti-linear involution of
$V_L.$
\begin{theorem}\label{latt}
Let $L$ be a positive definite integral lattice and $\phi$ be the anti-linear involution of $V_L$ defined above. Then $(V_L, \phi)$ is a unitary vertex operator superalgebra.
\end{theorem}
\pf  Since  the lattice vertex operator superalgebra $V_L$ is generated by
$$\{\alpha(-1) | \alpha\in L \}\cup \{e^\alpha |\alpha\in
L\},$$
by Proposition \ref{basic} it is sufficient to prove
the following identities
\begin{equation}\label{1.1}
(Y(e^{zL(1)}(-1)^{L(0)+2L(0)^2}z^{-2L(0)}\alpha(-1)\cdot1, z^{-1})w_1,
w_2)=(w_1, Y(\phi(\alpha(-1)\cdot 1), z)w_2),\end{equation}
 \begin{equation}\label{1.2}
 (Y(e^{zL(1)}(-1)^{L(0)+2L(0)^2}z^{-2L(0)}e^{\alpha},
z^{-1})w_1, w_2)=(w_1, Y(\phi(e^{\alpha}), z)w_2),\end{equation}
hold for any $w_1,w_2 \in V_L$.

 Assume that
$w_1= u\otimes
e^{\gamma_1},$ $w_2= v\otimes e^{\gamma_2}$ for some $u, v\in M(1)$ and $\gamma_1, \gamma_2\in L $.  By the
definition of the Hermitian form, we have
\begin{align*}
&(Y(e^{zL(1)}(-1)^{L(0)+2L(0)^2}z^{-2L(0)}\alpha(-1), z^{-1})w_1, w_2)\\
&=-z^{-2}\sum_{n\in \Z}(\alpha(n)w_1, w_2)z^{n+1}\\
&=\sum_{n\in \Z}-(w_1, \alpha(-n)w_2)z^{n-1}\\
&=(w_1, Y(\phi(\alpha(-1)), z)w_2),
\end{align*}
as desired.

To prove the identity (\ref{1.2}), note that we only need to consider the case that $\alpha+\gamma_1= \gamma_2$, by Lemma \ref{lattv},
\begin{align*}
&(Y(e^{zL(1)}(-1)^{L(0)+2L(0)^2}z^{-2L(0)}e^{\alpha}, z^{-1})w_1, w_2)\\
&=(-1)^{\frac{(\alpha, \alpha)+(\alpha, \alpha)^2}{2}}z^{-(\alpha, \alpha)}(Y(e^{\alpha}, z^{-1}) u\otimes
e^{\gamma_1}, v\otimes e^{\gamma_2})\\
&=(-1)^{\frac{(\alpha, \alpha)+(\alpha, \alpha)^2}{2}}z^{-(\alpha, \alpha)}(E^-(-\alpha,
z^{-1})E^+(-\alpha, z^{-1})e_{\alpha}(z^{-1})^{\alpha}u\otimes
e^{\gamma_1}, v\otimes
e^{\gamma_2})\\
&=z^{-(\alpha, \alpha)}(u\otimes
e^{\gamma_1},E^-(\alpha, z)E^+(\alpha,
z)(z^{-1})^{\alpha}e_{-\alpha}v\otimes e^{\gamma_2})\\
&=z^{-(\alpha, \alpha)}(u\otimes
e^{\gamma_1},E^-(\alpha, z)E^+(\alpha,
z)e_{-\alpha}(z^{-1})^{\alpha}(z^{-1})^{(\alpha,-\alpha)}v\otimes e^{\gamma_2})\\
&=(u\otimes e^{\gamma_1},E^-(\alpha,
z)E^+(\alpha,
z)e_{-\alpha}(z^{-1})^{\alpha}v\otimes e^{\gamma_2})\\
&=(w_1, Y(e^{-\alpha}, z)w_2)\\
&=(w_1, Y(\phi(e^{\alpha}), z)w_2).
\end{align*}
Hence, $(V_L, \phi)$ is a unitary vertex operator superalgebra.\qed
\section{Structures of unitary vertex operator superalgebras}
\def\theequation{3.\arabic{equation}}
\setcounter{equation}{0}
\subsection{Unitary structures of the direct sums of  unitary vertex operator superalgebras}
In this subsection, we shall show that there exists a unitary vertex operator superalgebra structure on the direct sum of unitary vertex operator superalgebras. We first recall from \cite{FHL}, \cite{LL} some facts about the direct sum of vertex operator superalgebras. Let $(V_1,\1_1, \w_1, Y_1), \cdots, (V_r,\1_r, \w_r, Y_r)$ be vertex operator superalgebras of the same central charge. Consider the vector space $V=V_1\oplus\cdots\oplus V_r$ and define the linear map $Y(\cdot, z)$ from $V$ to $(\End~ V)[[z, z^{-1}]]$ by
\begin{align*}
Y((v^{(1)},\cdots, v^{(r)}), z)=(Y_1(v^{(1)}, z), \cdots, Y_r(v^{(r)}, z))
\end{align*}
for $v^{(i)}\in V_i$, $1\leq i\leq r$. Set $\1=(\1_1, \cdots, \1_r)$ and $\w=(\w_1, \cdots, \w_r)$. It was proved in \cite{FHL}, \cite{LL} that $(V, \1, \w, Y)$ is a vertex operator superalgebra, the so-called direct sum vertex operator superalgebra of $(V_1,\1_1, \w_1, Y_1), \cdots, (V_r,\1_r, \w_r, Y_r)$. 
\begin{theorem}\label{direct}
Let $(V_1,\1_1, \w_1, Y_1), \cdots, (V_r,\1_r, \w_r, Y_r)$ be  vertex operator superalgebras of the same central charge and $V$ be the direct sum of $V_i$, $1\leq i\leq r$. Then $V$ has a unitary vertex operator superalgebra if each $V_i$, $1\leq i\leq r$, has a unitary vertex operator superalgebra structure.
\end{theorem}
\pf Assume that $(V_1,\1_1, \w_1, Y_1), \cdots, (V_r,\1_r, \w_r, Y_r)$ are unitary vertex operator superalgebras of the same central charge. Let $(,)_1, \cdots, (,)_r$ be the Hermitian forms on $V_1, \cdots, V_r$, respectively, and $\phi_1, \cdots, \phi_r$ be the anti-linear involutions of $V_1, \cdots, V_r$, respectively. We define a Hermitian form $(,)_V$ and an anti-linear involution $\phi$ of $V$ as follows:
\begin{align*}
&((u^{(1)},\cdots, u^{(r)}),(v^{(1)},\cdots, v^{(r)}))_V=\sum_{i=1}^r(u^{(i)}, v^{(i)})_{i},\\
&\phi((u^{(1)},\cdots, u^{(r)}))=(\phi_1(u^{(1)}),\cdots, \phi_r(u^{(r)})),
\end{align*}
for any $(u^{(1)},\cdots, u^{(r)}),(v^{(1)},\cdots, v^{(r)})\in V$. It is clear that $(,)_V$ is a positive definite Hermitian form on $V$ and that $\phi$ is an anti-linear involution of $V$. 

Let $L(n), L_i(n), ~1\leq i\leq r,$ be the operators defined by \begin{align*}
&Y_V(\w,z)=\sum_{n\in\Z}L(n)z^{-n-2},\\
&Y_i(\w_i,z)=\sum_{n\in\Z}L_i(n)z^{-n-2},~1\leq i\leq r,
 \end{align*}respectively. We then need to verify that the invariant property:
 \begin{align*}
 &(Y(e^{zL(1)}(-1)^{L(0)+2L(0)^2}z^{-2L(0)}(u^{(1)},\cdots, u^{(r)}), z^{-1})(v^{(1)},\cdots, v^{(r)}), (w^{(1)},\cdots, w^{(r)}))_V\\
 &\ \ =((v^{(1)},\cdots, v^{(r)}), Y(\phi((u^{(1)},\cdots, u^{(r)})), z)(w^{(1)},\cdots, w^{(r)}))_V,
\end{align*}
holds for any $(u^{(1)},\cdots, u^{(r)}),(v^{(1)},\cdots, v^{(r)}), (w^{(1)},\cdots, w^{(r)})\in V$. By the definitions, 
\begin{align}\label{udirect}
&(Y(e^{zL(1)}(-1)^{L(0)+2L(0)^2}z^{-2L(0)}(u^{(1)},\cdots, u^{(r)}), z^{-1})(v^{(1)},\cdots, v^{(r)}), (w^{(1)},\cdots, w^{(r)}))_V \nonumber\\
 &\ \ =(Y((e^{zL_1(1)}(-1)^{L_1(0)+2L_1(0)^2}z^{-2L_1(0)}u^{(1)},\cdots, e^{zL_r(1)}(-1)^{L_r(0)+2L_r(0)^2}z^{-2L_r(0)}u^{(r)}), z^{-1}) \nonumber \\
 &\ \ \ \ \ \ \cdot(v^{(1)},\cdots, v^{(r)}), (w^{(1)},\cdots, w^{(r)}))_V \nonumber \\
 &\ \ =\sum_{i=1}^r(Y_i(e^{zL_i(1)}(-1)^{L_i(0)+2L_i(0)^2}z^{-2L_i(0)}u^{(i)}, z^{-1})v^{(i)}, w^{(i)})_i\nonumber \\
 &\ \ =\sum_{i=1}^r(v^{(i)}, Y_i(\phi_i(u^{(i)}), z)w^{(i)})_i \nonumber \\
 &\ \ =((v^{(1)},\cdots, v^{(r)}), Y((\phi_1(u^{(1)}),\cdots, \phi_r(u^{(r)})), z)(w^{(1)},\cdots, w^{(r)}))_V \nonumber \\
 &\ \ =((v^{(1)},\cdots, v^{(r)}), Y(\phi((u^{(1)},\cdots, u^{(r)})), z)(w^{(1)},\cdots, w^{(r)}))_V,
\end{align}
as desired.
\qed
\subsection{Structures of unitary vertex operator superalgebras}
 Recall from \cite{DH} that a vertex operator superalgebra $V$ is called a {\em strong CFT type} vertex operator superalgebra if the following conditions hold:\\
(1)  $V_0=\C\1$ and $V_n=0$ if $n< 0$;\\
(2) $L(1)V_1=0$.

In this subsection, we shall show that any unitary vertex operator superalgebra is a direct sum of  some strong CFT type unitary simple vertex operator superalgebras with the same central charge.
\begin{theorem}\label{CFT}
Let $V$ be a unitary vertex operator superalgebra. Then there exist strong CFT type unitary simple vertex operator superalgebras $V^1, \cdots, V^k$ such that $V$ is isomorphic to the direct sum of the vertex operator superalgebras $V^1, \cdots, V^k$.
\end{theorem}
\pf First, note that $V$ is a unitary module for the Virasoro algebra, it follows immediately that $V_n=0$ for $n< 0$.

We next prove that if $V$ is a unitary vertex operator superalgebra, then $V$ is a direct sum of  some vertex operator superalgebras with the same central charge. Since $V$ is a unitary vertex operator superalgebra, we know that $V$ viewed as a $V$-module is completely reducible, then there exist irreducible $V$-submodules $V^1,\cdots, V^k$ of $V$ such that $V\cong V^1\oplus \cdots\oplus V^k$. 
Moreover, we have $(V^i, V^j)=0$ if $i\neq j$.
Consider the weight one subspaces $V^i_0$, $1\leq i\leq k$, of $V^i$, $1\leq i\leq k$, respectively. Since $V$ viewed as a $V$-module is generated by the vacuum vector, we know that each $V^i_0$, $1\leq i\leq k$, is not zero. We now prove that each $V^i_0$, $1\leq i\leq k$, is one-dimensional. Note that for any $u\in V^i_0$, we have
\begin{align*}
(L(-1)u, L(-1)u)=(u, L(1)L(-1)u)=(u, 2L(0)u)=0, 
\end{align*}
implies $L(-1)u=0$. Hence, $0=Y(L(-1)u, z)=\frac{d}{d z}Y(u, z)$, forces that $Y(u, z)=u_{-1}$. This further implies that $Y(v, z)u=e^{L(-1)z}Y(u, -z)v=e^{L(-1)z}u_{-1}v,$ for any $v\in V$. As a result, $v_ju=0$ for any $v\in V$ and $j\in \Z_{\geq 0}$. Thus, $[v_n, u_{-1}]=\sum_{j\geq 0} (v_ju)_{n-1-j}=0$. Since $V$ has countable dimension and $V^i$ is an irreducible $V$-module, we have $u_{-1}$ acting on $V^i$ is equal to $\lambda_u \id_{V^i}$ for some complex number $\lambda_u$ by the Schur's Lemma (cf. \cite{DY}). 
As a result, for any $u, v\in V^i_0$,  we have $\lambda_v u=v_{-1}u=u_{-1}v=\lambda_u v$, implies that $V^i_0$ is one-dimensional. Let $e^i\in V^i$,  $1\leq i\leq k$, be the vectors such that $e^i$ acting on $V^i$ is equal to $\id_{V^i}$. Then we  know that $(V^i, e^i, Y)$,  $1\leq i\leq k$, are  vertex superalgebras.

For any $n\in \frac{1}{2}\Z$, let $V^i_n=V^i_n\cap V$, $1\leq i\leq k$. Then we know that $V^i=\oplus_{n\in \frac{1}{2}\Z} V^i_n$, $1\leq i\leq k$. We further let $\w^i=L(-2)e^i$, $1\leq i\leq k$. Consider the operator $L^i(n)=\w^i_{n+1}$, we have $L^i(n)=\w^i_{n+1}=(L(-2)e^i)_{n+1}=(\w_{-1}e^i)_{n+1}=\w_{n+1}e^i_{-1}=L(n)e^i_{-1}$. As a result, the operators $L^i(n)$, $n\in \Z$, satisfy the Virasoro algebra relation with central charge $c=c_V$, where $c_V$ denotes the central charge of $V$. Moreover, $L^i(0)|_{V^i_n}=L(0)e^i_{-1}|_{V^i_n}=L(0)|_{V^i_n}=n\id_{V^i}$. Note also that for any $u, v\in V^i$,
\begin{align*}
Y(L^i(-1)u, z)v=Y(L(-1)e^i_{-1}u, z)v=Y(L(-1)u, z)v=\frac{d}{dz}Y(u, z)v.
\end{align*}Therefore, $(V^i, e^i, \w^i, Y)$,  $1\leq i\leq k$, are  vertex operator superalgebras. Since for any $v\in V_1$ and $u\in V_0$, we have $(L(1)v, u)=(v, L(-1)u)=0$, implies $L(1)v=0$. In particular, for any $v\in V^i_1$, we have $L^i(1)v=L(1)v=0$. Hence, $(V^i, e^i, \w^i, Y)$,  $1\leq i\leq k$, are strong CFT type vertex operator superalgebras. 

 Note that for any $1\leq i, j\leq k$ such that $i\neq j$, and any $v\in V^i, u\in V^j$, we have $e^{L(-1)z}Y(u, -z)v=Y(v, z)u\in V^i\cap V^j$. Hence, $Y(v, z)u=0$ for any $v\in V^i, u\in V^j$. This immediately implies that the vacuum vector $\1$ is equal to $e^1+\cdots+e^k$. Therefore, $\w=L(-2)\1=L(-2)(e^1+\cdots+e^k)=\w^1+\cdots+\w^k$. As a result, $V$ viewed as a vertex operator superalgebra is isomorphic to the direct sum of $V^i$, $1\leq i\leq k$.
 
 We proceed to prove that $(V^i, e^i, \w^i, Y)$,  $1\leq i\leq k$, are  unitary vertex operator superalgebras. Let $\phi$ be the anti-linear involution associated to the unitary structure. We first prove that $\phi(V^i)=V^i$, $1\leq i\leq k$. Note that for any $u, v\in V$,
\begin{align*}
(\phi(u), v)&=(\phi(u)_{-1}\1, v)\\
&=\Res_z z^{-1}(Y(\phi(u), z)\1, v)\\
&=\Res_zz^{-1}(\1, Y(e^{zL(1)}(-1)^{L(0)+2L(0)^2}z^{-2L(0)}u, z^{-1})v).
\end{align*}
Thus, $(\phi(V^i), V^j)=0$ if $i\neq j$. Note also that $(V^i, V^j)=0$ if $i\neq j$, which forces that $\phi(V^i)=V^i$, $1\leq i\leq k$. As a result, we get anti-linear involutions $\phi_i$, $1\leq i\leq k$, of $V^i$, $1\leq i\leq k$, respectively.  Now we can show that $(V^i, e^i, \w^i, Y)$,  $1\leq i\leq k$, are  unitary vertex operator superalgebras by the formula (\ref{udirect}).

 Finally, we prove that each  $V^i$, $1\leq i\leq k$,  is a simple vertex operator superalgebra. Otherwise, assume that $V^i$ has a nontrivial ideal $I$. Note that for any $u\in I, v\in V^i$,
  \begin{align*}
  (u, v)=\Res_z z^{-1}(u,Y(v, z)e^i )=\Res_z z^{-1} (Y(e^{zL^i(1)}(-1)^{L^i(0)+2L^i(0)^2}z^{-2L^i(0)}\phi (v), z^{-1})u, e^i).
  \end{align*}Since $(,)$ is positive definite, for any $u\in I$, there exists some $v\in V^i$ such that $(u, v)\neq 0$. This implies that $\Res_z z^{-1} (Y(e^{zL^i(1)}(-1)^{L^i(0)+2L^i(0)^2}z^{-2L^i(0)}\phi (v), z^{-1})u, e^i)\neq 0$. Set $w=\Res_z z^{-1} Y(e^{zL^i(1)}(-1)^{L^i(0)+2L^i(0)^2}z^{-2L^i(0)}\phi (v), z^{-1})u$, it follows from the discussion above that there exist some homogenous vectors $w_1, \cdots, w_n$ with different nonzero conformal weights such that $w=\lambda e^i+w_1+\cdots +w_n$ for some nonzero number $\lambda$. We now prove that $e^i\in I$, which is a contradiction.  Note that we have
   \begin{align*}
  & w=\lambda e^i+w_1+\cdots +w_n,\\
  &L^i(0)w=(\wt w_1) w_1+\cdots +(\wt w_n)w_n,\\
  &\ \ \ \ \ \ \ \ \ \ \ \ \ \vdots\\
  &L^i(0)^nw=(\wt w_1)^n w_1+\cdots +(\wt w_n)^nw_n,
   \end{align*}
   and $w, L^i(0)w, \cdots, L^i(0)^nw\in I$, it follows immediately from the Vandermende determinant formula  that $e^i\in I$.
   Hence, each  $V^i$, $1\leq i\leq k$, is a simple vertex operator superalgebra. The proof is complete.
   \qed
   \begin{corollary}
Let $(V, \phi)$ be a  unitary simple vertex operator superalgebra. Then the Hermitian form $(,)$ is uniquely determined by the value $(\1, \1)$.
\end{corollary}
\pf Since the unitary vertex operator superalgebra $V$ is simple, it follows from the argument in the proof of Theorem \ref{CFT} that $V$ is of CFT type. 
As a result, the Hermitian form is uniquely determined by the value $(\1, \1)$.
\qed

In the following, for any unitary simple vertex operator superalgebra $V$, we will
normalize the Hermitian form $(,)$ on $V$ such that $(\1, \1)=1$.
\section{Classification of unitary vertex operator superalgebras generated by the subspaces with small conformal weights}
\def\theequation{4.\arabic{equation}}
\setcounter{equation}{0}
\subsection{Classification of unitary vertex operator superalgebras generated by the subspaces with conformal weight $\frac{1}{2}$} In the last section, we have shown that the conformal weights of any unitary vertex operator superalgebra are nonnegative. Then it is natural to consider first the classification of unitary vertex operator superalgebras generated by the subspaces with conformal weight $\frac{1}{2}$. In subsection \ref{superhe}, we have proved that the vertex operator superalgebras associated to the Heisenberg superalgebras have unitary vertex operator superalgebra structures. Moreover, these vertex operator superalgebras are generated by the subspaces with conformal weight $\frac{1}{2}$.  In this subsection, we shall show that any unitary vertex operator superalgebras generated by the subspaces with conformal weight $\frac{1}{2}$ is a direct sum of  vertex operator superalgebras associated to  some  Heisenberg superalgebras. First,  we recall from \cite{X} a  characterization of  the vertex operator superalgebras associated to the Heisenberg superalgebras.
\begin{theorem}\label{char}
Let $V$ be a simple vertex operator superalgebra generated by the subspace $V_{\frac{1}{2}}$. 
Then $V_{\frac{1}{2}}$ has a nondegenerate symmetric bilinear form $\langle, \rangle$ such that $V$ is isomorphic to the vertex operator superalgebra associated to the Heisenberg superalgebra $\hat V_{\frac{1}{2}}$.
\end{theorem}

As a consequence, we have the following classification of unitary vertex operator superalgebras generated by the subspaces with conformal weight $\frac{1}{2}$.

\begin{theorem}
   Let $V$ be a unitary vertex operator superalgebra generated by the subspace $V_{\frac{1}{2}}$. 
   Then $V$ is isomorphic to a direct sum of  vertex operator superalgebras associated to some Heisenberg superalgebras.
   \end{theorem}
   \pf Since $V$ is a unitary vertex operator superalgebra, by Theorem \ref{CFT}, we know that there exists unitary simple vertex operator superalgebras $V^1, \cdots, V^k$ such that $V$ is the direct sum of $V^i$, $1\leq i\leq k$. By assumption, $V$ is generated by the subspace $V_{\frac{1}{2}}$ with conformal weight $\frac{1}{2}$, it follows that each $V^i$, $1\leq i\leq k$, is  generated by the subspace with conformal weight $\frac{1}{2}$. Therefore, by Theorem \ref{char}, we know that each $V^i$, $1\leq i\leq k$, is isomorphic to  vertex operator superalgebra associated to some Heisenberg superalgebra. The proof is complete.\qed
\subsection{The structures of the weight one subspaces of unitary vertex operator algebras} In the last subsection, we have classified the unitary vertex operator superalgebras generated by the subspaces with conformal weight $\frac{1}{2}$. Our next goal is to classify the unitary vertex operator algebras generated by the subspaces with conformal weight $1$. In this subsection, we shall study first the structures of the weight one subspaces of unitary vertex operator algebras. Let $V$ be a unitary vertex operator algebra such that $V_1\neq 0$. Note that for $u, v\in V_1$, we have $u_1v\in V_0$, it follows from Theorem \ref{CFT} that $L(-1)u_1v=0$. As a result, we know that $V_1$ has a Lie algebra structure with the Lie bracket defined by $[u, v]=u_0v$ for any $u, v\in V_1$ (cf. \cite{B}).  We next prove that $V_1$ is a reductive Lie algebra, that is, $V_1$ viewed as a $V_1$-module is completely reducible (cf. \cite{H}).
\begin{theorem}\label{reductive}
Let $V$ be a unitary vertex operator algebra. Then $(V_1, [,])$ is a reductive Lie algebra.
\end{theorem}
\pf  Let $(,)$ be the positive definite Hermitian form on $V$, then it follows from Theorem \ref{CFT} that $(L(1)v, u)=(v, L(-1)u)=0$ for any $v\in V_1$ and $u\in V_0$, which implies $L(1)v=0$. Thus, for any $u,v, w\in V_1$,
\begin{align}\label{invariant}
\nonumber([u, v], w)&=(u_0v, w)\\ \nonumber
&=\Res_z(Y(u, z)v, w)\\ \nonumber
&=\Res_z(v, Y(e^{zL(1)}(-z^{-2})^{L(0)}\phi(u), z^{-1})w)\\ \nonumber
&=\Res_z(v, Y((-z^{-2})\phi(u), z^{-1})w)\\
&=-(v, [\phi(u), w]).
\end{align}
Therefore, for any submodule $I$ of $V_1$, the orthogonal complement $I^{\perp}$ of $I$, which is defined by $I^{\perp}=\{u\in V_1|(u, v)=0\text{ for any }v\in I\}$, is also a submodule of $V_1$. This implies that $V_1$ viewed as a $V_1$-module is completely reducible. This completes the proof.\qed
\begin{remark}
The similar result in Theorem \ref{reductive} has been established in \cite{DM1} for the vertex operator algebra $V$ satisfying the following conditions:\\
(1) $V$ is simple, rational and $C_2$-cofinite;\\
(2) $V$ is of strong CFT type.
\end{remark}

\subsection{A characterization of the unitary vertex operator algebra associated to the affine Lie algebra}
We now assume that $V$ is a unitary simple vertex operator algebra such that $V_1$ contains a simple Lie algebra $\g$. In the following, we shall study the structure of the vertex subalgebra $U$ generated by $\g$.
 
 By assumption, we know that $V$ is a self-dual,  strong CFT, simple vertex operator algebra. Then there is a unique nondegenerate invariant bilinear form $\langle ,\rangle$ on $V$ such that $\langle\1,\1\rangle=-1$ (see \cite{Li1}). Moreover, the nondegenerate invariant bilinear form $\langle ,\rangle$ is symmetric (see \cite{FHL}). Note that we have defined in Proposition \ref{self-dual} a nondegenerate invariant bilinear form $\langle ,\rangle$ on $V$ such that $\langle u,v\rangle=-(u, \phi(v))$ for any $u,v\in V$ and that $\langle \1,\1\rangle=-(\1, \1)=-1$. Therefore, the nondegenerate invariant bilinear form $\langle ,\rangle$ defined by  $\langle u,v\rangle=-(u, \phi(v))$ is symmetric.
 
 By restricting the bilinear form $\langle ,\rangle$ to $V_1$, we then get a nondegenerate symmetric bilinear on $V_1$. We proceed to prove $\langle,\rangle$ is an invariant form, that is, $\langle [u, v], w\rangle=\langle u, [v, w]\rangle$ for any $u, v, w\in V_1$. By the formula (\ref{invariant}), we have $\langle [u, v], w\rangle=-([u, v], \phi(w))=-(u, [\phi(v), \phi(w)])=-(u, \phi([v, w]))=\langle u, [v, w]\rangle$. Hence, $\langle,\rangle$ is a nondegenerate invariant symmetric bilinear form on $V_1$. Since $V_1$ is reductive, the restriction of $\langle,\rangle$ on $\g$ is still nondegenerate.  
 This implies that there exists a nonzero complex number $c$ such that $\langle,\rangle|_{\g\times \g}=c\kappa (,)$, where $\kappa(,)$ denotes the normalized Killing form of $\g$ such that  $\kappa(\theta, \theta)=2$ for the highest root $\theta$ of $\g$.
 
 Note also that for any $u,v\in V_1$, \begin{align*}
(u, v)&=(u, v_{-1}\1)=\Res_zz^{-1}(u, Y(v, z)\1)\\
&=\Res_zz^{-1}(Y(e^{zL(1)}(-z^{-2})^{L(0)}\phi(v), z^{-1})u, \1)\\
&=\Res_z-z^{-3}(Y(\phi(v), z^{-1})u, \1)\\
&=-(\phi(v)_1u, \1).
\end{align*} It follows that $u_1v=-(v, \phi(u))\1=\langle v, u\rangle\1=\langle u, v\rangle\1$. Thus, we have
\begin{align*}
[u_m, v_n]&=\sum_{i\geq 0}{m\choose i}(u_iv)_{m+n-i}\\
&=(u_0v)_{m+n}+m(u_1v)_{m+n-1}\\
&=[u, v]_{m+n}+m\langle u, v\rangle\delta_{m+n, 0}\id.
\end{align*}
Since $U$ is generated by $\g$, we have in fact proved that $U$ is a highest weight module for the affine Lie algebra $\hat{\g}$ with highest weight vector $\1$ and level $c$. In particular, $U$ viewed as a $\hat{\g}$-module is isomorphic to some quotient module of the Verma module $M_{\hat{\g}}(c, 0)$. Furthermore, we have the following
\begin{theorem}\label{affine}
Let $U$ be the vertex  subalgebra defined above. Then $U$ viewed as a $\hat{\g}$-module is irreducible, i.e., $U$ isomorphic to  $L_{\hat{\g}}(c, 0)$ for some complex number $c$.
\end{theorem}
\pf By assumption, $V$ has a positive definite Hermitian form such that
$$(Y(e^{zL(1)}(-z^{-2})^{L(0)}a, z^{-1})u, v)=(u,
Y(\phi(a), z)v),$$holds for any $a, u, v\in V$. In particular, we have for any $a\in \g$, $u, v\in U$ and $n\in \Z$,
$$(a_nu, v)=-(u,
\phi(a)_{-n}v).$$ Thus we can prove that the radical of $(,)$ is a $\hat{\g}$-module and that the radical of $(,)$ contains the maximal proper submodule which has zero intersection with highest weight subspace. Since the radical of $(,)$ is equal to zero, we get that the maximal proper submodule which has zero intersection with highest weight subspace is zero. Therefore, $U$ is irreducible.\qed
\vskip0.25cm
We proceed to prove that $U$ viewed as a $\hat{\g}$-module is an integrable module, that is, $c$ is a positive integer (cf. \cite{K}).
\begin{theorem}\label{simplelie}
Let $U$ be the vertex subalgebra defined above. Then $U$ viewed as a $\hat{\g}$-module is integrable. 
\end{theorem}
\pf Let $\phi$ be the anti-linear involution of $V$ associated to the unitary structure. We first prove that $\phi(\g)=\g$. Note that we have proved in Theorem \ref{reductive} that $V_1$ is reductive and then $V_1\cong \g_1\oplus \cdots\oplus \g_k\oplus \h$, where $\g_i$, $1\leq i\leq k$, are simple Lie subalgebras of $V_1$ and $\h$ is an abelian Lie subalgebra of $V_1$. Moreover, we have $\g_i$, $1\leq i\leq k$, and $\h$ are mutually orthogonal with respect to $(,)$. On the other hand, we have proved that the bilinear form $\langle ,\rangle$ defined by  $\langle u,v\rangle=-(u, \phi(v))$ is  a nondegenerate invariant symmetric bilinear form. Since $\g_i$, $1\leq i\leq k$, are simple Lie algebras, it follows that $\g_i$, $1\leq i\leq k$, and $\h$ are mutually orthogonal with respect to $\langle ,\rangle$. As a result, $\phi(\g_i)=\g_i$, $1\leq i\leq k$, and $\phi(\h)=\h$.

Note also that $\phi([u, v])=\phi(u_0v)=\phi(u)_0\phi(v)=[\phi(u), \phi(v)]$ for any $u, v\in V_1$. Thus, $\phi$ induces an anti-linear involution of $\g$, which is still denoted by $\phi$. Hence, $\phi$ is a real structure of $\g$ (cf. \cite{O}). Consider the real form $\g^\phi=\{v\in \g|\phi(v)=v\}$ of $\g$ with respect to $\phi$. We then prove that $\g^\phi$ is a compact Lie algebra, that is, there exists a positive definite symmetric bilinear form on $\g^\phi$ (cf. \cite{O}). Note that we have proved that the bilinear form $\langle, \rangle$ defined by $\langle u, v\rangle=(u, \phi(v))$, for any $u, v\in \g$, is a symmetric invariant bilinear form, 
we then get a  positive definite symmetric bilinear form on $\g^\phi$. Therefore, $\g^\phi$ is a compact Lie algebra and $\phi$ is a compact real structure of $\g$. Since the compact real structure of a simple Lie algebra is unique up to conjugacy by inner automorphism, we can choose a Cartan subalgebra $\h_1$ of $\g$ such that $\phi$ is determined by $\phi(e_i)=-f_i, \phi(f_i)=-e_i, \phi(h)=-h, h\in \h_1$, where $e_i, f_i$ denote the Chevalley generators of $\g$ with respect to $\h_1$ (cf. \cite{O}).

We then consider the compact involution $\hat \phi$ of $\hat{\g}$, which is the anti-linear map $\hat \phi: \hat{\g}\to \hat{\g}$ determined by $\hat\phi(u\otimes t^n)=\phi(u)\otimes t^{-n}, u\in \g,$ and $\hat \phi(C)=-C$ (cf. \cite{K}).
Note that with respect to $\hat \phi$ and the positive definite Hermitian form $(,)$, the $\hat{\g}$-module $U$ is a unitarizable module, i.e., 
 $(x\cdot u, v)=-(u, \hat{\phi}(x)\cdot v)$ for any $x\in \hat{\g}$, $u, v\in U$ (cf. \cite{K}). As a result, $U$ viewed as $\hat{\g}$-module is integrable. The proof is complete. \qed
 
 We now in a position to give a characterization of the unitary vertex operator algebra associated to the affine Lie algebra.
 \begin{theorem}\label{intvoa}
 Let $V$ be a unitary vertex operator algebra such that $V_1$ is a simple Lie algebra and that $V$ is generated by $V_1$. Then $V$ is isomorphic to the affine vertex operator algebra $L_{\hat{V_1}}(c, 0)$.
 \end{theorem}
 \pf We first prove that $V$ is a simple vertex operator algebra. Since $V$ is a unitary vertex operator algebra, by Theorem \ref{CFT}, we know that there exist unitary simple vertex operator algebras $V^1, \cdots, V^k$ such that $V$ is the direct sum of $V^i$, $1\leq i\leq k$. However, by assumptions, $V$ is generated by the subspace $V_{1}$ and $V_1$ is a simple Lie algebra, it follows that $V$ is isomorphic to some $V^j$, $1\leq j\leq k$. Therefore, by Theorem \ref{CFT}, we know that $V$ is a simple vertex operator algebra.
 
 We now only need to prove that the conformal vector $\w$ of $V$ is equal to that of $L_{\hat{V_1}}(c, 0)$. Consider the component operators $L(n)=\w_{n+1}$, $n\in \Z$, we have: For any $x\in V_1$, $n, m\in \Z$,
 \begin{align*}
 [L(n), x_m]&=[\w_{n+1}, x_m]\\
 &=\sum_{i\geq 0} {n+1\choose i}(\w_ix)_{m+n+1-i}\\
 &=(\w_0x)_{m+n+1}+(n+1)(\w_1x)_{m+n}\\
 &=-mx_{m+n}.
 \end{align*}
 On the other hand, let $x_1,\cdots,  x_d$ be an orthonormal basis of $V_1$ with respect to $\kappa(, )$, since $c$ is a positive integer, we can define $\w'=\frac{1}{2(c+h)}\sum_{i=1}^dx_i(-1)^2\1$, where $h$ denotes the dual Coxeter number of $V_1$. It is known that the component operators $L'(n)=\w'_{n+1}$, $n\in \Z$, satisfy $[L'(n), x_m]=-mx_{m+n}$ for any $x\in V_1$, $n, m\in \Z$ (see \cite{LL}). As a consequence, we have $\w=\w'$ by Theorem 6.3 in \cite{Lian}. The proof is complete.\qed
\begin{remark}
The structure of the vertex operator algebra $V$ satisfying the following conditions:\\
(1) $V$ is simple and $C_2$-cofinite,\\
(2) $V$ is of strong CFT type,\\
(3) $V_1$ contains a simple Lie subalgebra,\\
has been studied in \cite{DM2} and 
the similar results in Theorems \ref{affine}, \ref{simplelie},  \ref{intvoa} have been established for the vertex operator algebra $V$ satisfying these conditions.
\end{remark}
\subsection{Classification of unitary vertex operator algebras generated by the weight one subspaces} In the last subsection, we have studied the structure of the unitary simple vertex operator algebra $V$ under the condition that $V_1$ contains a simple Lie subalgebra. To classify unitary vertex operator algebras generated by the weight one subspaces, it remains to study the structure of the unitary simple vertex operator algebra $V$ under the condition that $V_1$ contains an abelian Lie algebra. In this subsection, we shall study this case and then classify unitary vertex operator algebras generated by the weight one subspaces.

Let $V$ be a unitary simple vertex operator algebra such that $V_1\cong \g_1\oplus \cdots\oplus \g_k\oplus \h$, where $\g_i$, $1\leq i\leq k$, are simple Lie subalgebras of $V_1$ and $\h$ is an abelian Lie subalgebra of $V_1$. Our goal is to study the structure of the vertex subalgebra $U$ generated by $\h$. Note that we have proved in Theorem \ref{simplelie} that there is a nondegenerate symmetric bilinear form $\langle , \rangle$ on $\h$ defined by $\langle u, v\rangle=-(u, \phi (v))$ for any $u, v\in \h$. Then we can consider the Heisenberg algebra $\hat{\h}=\h\otimes \C[t, t^{-1}]\oplus \C K$ with the communication relations: For any $u, v\in \h$,
\begin{align*}
&[u(m), v(n)]=m\langle u, v\rangle\delta_{m+n, 0}K,\\
&[K, x]=0, \text{ for any } x\in \hat{\h}.
\end{align*}
Since $U$ is generated by $\h$, we can prove that $U$ is a highest weight module for the Heisenberg Lie algebra $\hat{\h}$ with highest weight vector $\1$ and level $1$. In particular, $U$ viewed as an $\hat{\h}$ module is isomorphic to the irreducible highest weight $\hat{\h}$-module $V_{\hat{V_1}}(1, 0)$ (cf. \cite{FLM}, \cite{LL}). Moreover, we have
\begin{theorem}\label{abelian}
Let $V$ be a unitary simple vertex operator algebra such that $V$ is generated by $V_1$ and that $V_1$  is an abelian algebra. Then the vertex operator algebra $V$ is isomorphic to $V_{\hat{V_1}}(1, 0)$, where we use $V_{\hat{\h}}(1, 0)$ to denotes the Heisenberg vertex operator algebra associated to a Heisenberg algebra $\hat{\h}$.
\end{theorem}
\pf It is good enough to prove that the conformal vector $\w$ of $V$ coincides with that of $V_{\hat\h}(1,0)$. First, note that we have proved in Theorem \ref{intvoa} that $[L(n), x_m]=-mx_{m+n}$ holds for any $x\in V_1$, $n, m\in \Z$, where $L(n)=\w_{n+1}$. On the other hand, let $\alpha_1, \cdots, \alpha_d$ be an orthonormal basis of $V_1$ with respect to $\langle, \rangle$ and consider $\w'=\frac{1}{2}\sum_{i=1}^d\alpha_i^2(-1)\1$. It is known that $[L'(n), x_m]=-mx_{m+n}$ holds for any $x\in V_1$, $n, m\in \Z$, where $L'(n)=\w'_{n+1}$ (cf. \cite{LL}).  As a consequence, we have $\w=\w'$ by Theorem 6.3 in \cite{Lian}. The proof is complete.\qed

\vskip.25cm
We now in a position to classify the unitary vertex operator algebras generated by the weight one subspaces. Let $V$ be a unitary vertex operator algebra such that $V$ is generated by $V_1$. By Theorem \ref{CFT}, we know that $V$ is a direct sum of unitary simple vertex operator algebras $V^1, \cdots, V^k$. Since $V$ is generated by $V_1$, we know that the weight one subspaces $V^i_1$, $1\leq i\leq k$, of $V^i$, $1\leq i\leq k$, respectively, are nonzero and that each $V^i$ is generated by $V^i_1$. As a consequence, we know that $V^i_1$, $1\leq i\leq k$, are reductive Lie algebras. Assume that $V^i_1$, $1\leq i\leq k$, are isomorphic to $\g^i_1\oplus\cdots\oplus \g^i_{l_i}\oplus \h^i$, $1\leq i\leq k$, respectively, where $\g^i_j$ , $1\leq i\leq k$, $1\leq j\leq l_i$, are simple Lie algebras and $\h^i$, $1\leq i\leq k$, are abelian Lie algebras. Then we have
\begin{theorem}\label{all}
The unitary vertex operator algebra $V$ is isomorphic to 
\begin{align*}
L_{\hat{\g^1_1}}(c^1_1, 0)\otimes\cdots\otimes L_{\hat{\g^1_{l_1}}}(c^1_{l_1}, 0)\otimes V_{\hat\h^1}(1, 0)\oplus \cdots\oplus L_{\hat{\g^k_1}}(c^k_1, 0)\otimes\cdots\otimes L_{\hat{\g^k_{l_k}}}(c^k_{l_k}, 0)\otimes V_{\hat\h^k}(1, 0)
 \end{align*} for some positive integers $c^i_j$, $1\leq i\leq k$, $1\leq j\leq l_i$.
\end{theorem}
\pf We now only need to show that the conformal vectors of $V^i$, $1\leq i\leq k$, coincide with those of $L_{\hat{\g^i_1}}(c^i_1, 0)\otimes\cdots\otimes L_{\hat{\g^i_{l_1}}}(c^i_{l_i}, 0)\otimes V_{\hat\h^i}(1, 0)$, $1\leq i\leq k$, respectively.  This follows immediately by the same arguments in the proofs of Theorems \ref{intvoa}, \ref{abelian}. The proof is complete.\qed
\begin{remark}
The similar results in Theorems \ref{abelian}, \ref{all} have been established in \cite{DM2} for the vertex operator algebra $V$ satisfying the following conditions:\\
(1) $V$ is simple and $C_2$-cofinite;\\
(2) $V$ is of strong CFT type.
\end{remark}

\subsection{A characterization of the lattice vertex operator algebra}In this subsection, we use the results obtained in the last subsection to give a characterization of the lattice vertex operator algebra. We will use the notations in the last subsection.
\begin{theorem}
Let $V$ be a  unitary simple vertex operator algebra such that $V_1$ is an abelian Lie algebra and $\dim V_1=c$, where $c$ denotes the central charge of $V$. Then there exist an abelian subalgebra $\h_1$ of $V_1$ and a  positive definite even lattice $L$ in $V_1$ such that $V$ is isomorphic to $V_{\hat\h_1}(1,0)\otimes V_{L}$.
\end{theorem}
\pf By the discussion above, we know that $V$ is a module for the Heisenberg algebra $\hat V_1$. Thus, by the Stone-Von Neumann Theorem we know that $V$ viewed as a $\hat V_1$-module has the decomposition
$V=V_{\hat V_1}(1, 0)\otimes \Omega_V,$
where $V_{\hat V_1}(1, 0)$ is the highest weight $\hat V_1$-module with highest weight $(1, 0)$ and $\Omega_V=\{v\in V|x(m)v=0, \text{ for any } x\in V_1,\ m>0\}$.

We next show that $x(0)$ acts semisimplely on $V$ for any $x\in V_1$. Note that the action on $V_n$ is closed for any $n\in \Z$, so it is good enough to prove that $x(0)$ acts semisimplely on $V_n$ for any $x\in V_1$. Let $(,)$ be the positive definite Hermitian form on $V$ and $\phi$ be the anti-linear involution of $V$ associated to the unitary structure. Then we have $(x(0)u, v)=(u, -\phi(x)(0)v)$ for any $u, v\in V_n$. Since $\phi$ is an anti-linear involution, by viewing $V_1$ as a vector space $V_{1,\R}$ over the field $\R$ of real numbers, we can view $\phi$ as a linear involution of $V_{1,\R}$. Thus, $V_{1,\R}$ has a decomposition $V_{1,\R}=V_{1,\R}^+\oplus V_{1,\R}^-$, where $V_{1,\R}^+=\{v\in V_{1,\R}| \phi(v)=v\}$ and $V_{1,\R}^-=\{v\in V_{1,\R}| \phi(v)=-v\}$. 
We first assume that $x$ is an element in $V_{1,\R}^-$, then we have $(x(0)u, v)=(u, x(0)v)$ for any $u, v\in V_n$. It follows that $x(0)$ acts semisimplely on $V_n$. We next assume that $x\in V_{1,\R}^+$, then we have $\phi(ix)=-i\phi(x)=-ix$, that is, $ix\in V_{1,\R}^-$. Therefore, $ix(0)$ acts semisimplely on $V_n$. Thus, for any $x\in V_1$, we have $x(0)$ acts semisimplely on $V_n$.

As a result, for any $x\in V_1$, $x(0)$ acts semisimplely on $\Omega_V$. For any $\lambda\in V_1$, let $\Omega_V(\lambda)=\{v\in \Omega_V|x(0)v=\langle x, \lambda\rangle v\text{ for any } x\in V_1\}$, then we have
$\Omega_V=\oplus_{\lambda\in V_1}\Omega_V(\lambda).$
Let $L=\{\lambda\in V_1|\Omega_V(\lambda)\neq 0\}$. Since $V$ is a unitary simple vertex operator algebra, by a similar argument in the proof of Theorem 2 in \cite{DM1} we can prove that $L$ is an additive group. We proceed to prove that $L$ is a positive definite even lattice. Let $\alpha_1, \cdots, \alpha_c$ be an orthonormal basis of $V_1$ with respect to $\langle, \rangle$. Then we know that $\w'=\frac{1}{2}\sum_{i=1}^c\alpha_i^2(-1)\1$ is the conformal vector of the Heisenberg vertex operator algebra $V_{\hat V_1}(1, 0)$. Set $\w''=\w-\w'$, where  $\w$ is the conformal vector of $V$. We now prove that $\w''$ is a Virasoro element. Set $L(n)=\w_{n+1}$ and $L'(n)=\w'_{n+1}$, $n\in \Z$. Then we have $[L(1), x_{-1}]=x_0$ for any $x\in V_1$. Therefore, $L(1)\w'=0$. It follows that $\w''$ is a Virasoro element with central charge $0$ (see \cite{FZ}). Set
$L''(n)=\w''_{n+1}$, then we have
 \begin{align*}
 (\w'', \w'')&=(\w-\w',\w-\w')\\
& =((L(-2)-L'(-2))\1, (L(-2)-L'(-2))\1)\\
& =(\1, (L(2)-L'(2))(L(-2)-L'(-2))\1)\\
& =(\1, L''(2)L''(-2)\1)=0.
 \end{align*}
 Since the Hermitian form is positive definite, $\w''=0$. As a consequence, $\w'$ is equal to the conformal vector of $V$. This forces that $\frac{\langle \lambda, \lambda\rangle}{2}\in \Z^+$ for any $\lambda\in L$ (see \cite{DM1}). Hence, $L$ is a positive definite even lattice.

We now set $\h_2=\C\otimes_{\Z}L$ and $\h_1=\{v\in V_1|\langle x, v\rangle=0 \text{ for any } x\in \h_2\}$. Then $V_1=\h_1\oplus \h_2$ and $\langle, \rangle|_{\h_1\times\h_1}$, $\langle, \rangle|_{\h_2\times\h_2}$ are nondegenerate. Thus we can consider the Heisenberg algebras $\hat \h_1$, $\hat \h_2$. Moreover, we have $V_{\hat V_1}(1, 0)=V_{\hat \h_1}(1, 0)\otimes V_{\hat \h_2}(1, 0)$. Thus, $V$ viewed as a module for $V_{\hat \h_1}(1, 0)\otimes V_{\hat \h_2}(1, 0)$ has the following decomposition:
\begin{align*}
V=\left(V_{\hat \h_1}(1, 0)\otimes V_{\hat \h_2}(1, 0)\right)\otimes \Omega_V=V_{\hat \h_1}(1, 0)\otimes\left(\oplus_{\lambda\in \h_2} V_{\hat \h_2}(1, 0)\otimes\Omega_V(\lambda)\right).
\end{align*}
It follows immediately from Proposition 5.3 of \cite{DM1} that the vertex operator algebra $V$ is isomorphic to $V_{\hat \h_1}(1, 0)\otimes V_L$. This completes the proof.\qed

As an immediate consequence, we have the following 
\begin{corollary}\label{clattice}
Let $V$ be a  unitary simple vertex operator algebra such that $V_1$ is an abelian Lie algebra and $\dim V_1=c$, where $c$ denotes the central charge of $V$. If we assume in addition that $V$ is rational or $C_2$-cofinite, then there exists a positive definite even lattice $L$ in $V_1$ such that $V$ is isomorphic to $V_L$.
 \end{corollary}
\begin{remark}
The similar result in Corollary \ref{clattice} has been established in \cite{DM1}, \cite{M} for the vertex operator algebra $V$ satisfying the following conditions:\\
(1) $V$ is simple, rational and $C_2$-cofinite;\\
(2) $V$ is of strong CFT type.\\
There is also a characterization of vertex algebras associated to even lattices in \cite{LX}.
\end{remark}
{\bf Acknowledgement}. The authors wish to thank Professor Ching Hung Lam for useful comments and
valuable suggestions.

\end{document}